\documentclass[12pt,thmsa]{article}
\usepackage{amsmath}
\usepackage{graphicx}
\usepackage{amsfonts}
\usepackage{amssymb}
\usepackage{mathrsfs}
\usepackage{dsfont}

\newtheorem{Assumption}{Assumption}

\newtheorem{lem}{Lemma}
\newtheorem{thm}{Theorem}

\newtheorem{remark}{Remark}

\begin{document}

\begin{center}
{\Large \textbf{Rates of strong uniform consistency for the $k$-nearest neighbors kernel estimators of density and regression function}}

\bigskip

Luran   BENGONO MINTOGO, Emmanuel de Dieu  NKOU  and  Guy Martial  NKIET 

\bigskip

\textsuperscript{}LPSI, URMI, Universit\'{e} des Sciences et Techniques de Masuku,  Franceville, Gabon.

\bigskip

E-mail adresses : luranbengono@gmail.com,   emmanueldedieunkou@gmail.com, guymartial.nkiet@univ-masuku.org.

\bigskip
\end{center}

\noindent\textbf{Abstract.}We adress the problem of consistency of the  $k$-nearest neighbors kernel  estimators of the density and the regression function in the multivariate case. We get the  rates of strong uniform consistency on the whole space $\mathbb{R}^p$ for these estimators under specified assumptions.

\bigskip


\section{Introduction}
Estimation  of the density and the  regression function are  important and classical issues  in nonparametric statistics which has been  intensively  addressed since many years, so leading to an abundant literature. Without a doubt, the most popular estimators that have been tackled in this context are the kernel estimators, namely the Parzen-Rosenblatt estimator of the density and the Nadaraya-Watson estimator of the regression function. However, the practical choice of the bandwith on which these estimators rely
is not straightforward and stills a challenging issue. This is why alternative estimators, which do not require to make  such a choice, have been proposed. Among them, the $k$-nearest neighbors ($k$-NN) kernel estimators have attracted particular attention. They have the same form than the kernel estimators, but with bandwith replaced by the Euclidean distance between the point to which the estimator is calculated and the $k$th nearest neighbor  of this point among the observations. Earlier works on this topic go back to \cite{moore77} for density estimation and to \cite{collomb80} for the case of regression function.   These papers established  the strong uniform consistency and the strong  pointwise consistency, respectively, of the tackled estimators.  After these works, some  others  studying various aspects related to the aforementioned  estimators  were introduced in the literature. For example, \cite{mack79} derived expressions describing the asymptotic behavior of the bias and variance of the $k$-NN density estimates, \cite{bhattacharya90} introduced an adaptative optimal choice of $k$ in multivariate $k$-NN density and regression estimation and  \cite{lu98} proved  strong pointwise consistency for the $k$-NN estimators of the density and the regression function in the context of $\alpha$-mixing stationnary sequences. The most recent works on the $k$-NN kernel estimators concern  the case of functional data (e.g., \cite{attouch17,burba08,kudraszow13,lian11}), and  that of spatial data (\cite{ahmed23}). Curiously, there is almost no work devoted to determining the convergence rates of the aforementioned $k$-NN kernel estimators in both the univariate and the multivariate cases. However, \cite{zhang90} derived  rates of strong uniform convergence, on any compact subset of $\mathbb{R}$,  of the $k$-NN kernel   density estimator, but only  in the univariate case.  To the best of our knowledge, there is no work dealing with derivation of such rates for the $k$-NN kernel estimator of the regression function  either in the univariate case or in the  multivariate case.

 In this paper, we address the problem of determining rates of strong uniform consistency for the $k$-NN kernel estimators  of multivariate density and regression function. In Section \ref{knnest}, we define the estimators that will be tackled. For the density, it is the usual $k$-NN kernel estimator but for the regression function, we slightly modify the classical one as it was done in \cite{zhu96}.  Section \ref{results} presents the used assumptions and gives the main results. The proofs of the theorems are postponed in Section  \ref{preuves}.

\section{The $k$-NN kernel estimators}\label{knnest}
Let $\left\{(X_i,Y_i)\right\}_{n\in\mathbb{N}^\ast}$ be an i.i.d. sample of a pair $(X,Y)$ of random variables valued into $\mathbb{R}^p\times\mathbb{R}$, with $p\geqslant 1$. We denote by $f$ the density of $X$ and, assuming that $\mathbb{E}\left(\left\vert Y\right\vert\right)<+\infty$, we consider the regression function $r$ defined as $$r(x)=\mathbb{E}\left(Y\vert X=x\right), \,\,\,x\in\mathbb{R}^p.$$ For estimating $f$, the $k$-NN kernel density estimator $\widehat{f}_n$  was introduced (see, e.g., \cite{moore77}); it is defined   as 
\begin{equation*}
\widehat{f}_{n}(x)=\frac{1}{n(R_{n}(x))^{p}}\sum\limits_{i=1}^{n}K\left(\frac{X_{i}-x}{R_{n}(x)}\right),\hspace{0.2cm}\textrm{$x\in  \mathbb{R}^{p}$},
\end{equation*}
where $K\,:\,\mathbb{R}^p\rightarrow\mathbb{R}$ is a multivariate kernel, and
\begin{equation*}
R_{n}(x)=\min\left\{h\in\mathbb{R}_{+}^{*}\bigg{/}\sum\limits_{i=1}^{n}\mathds{1}_{\mathcal{B}(x,h)}(X_{i})=k_{n}\right\}
\end{equation*}
with $\mathcal{B}(x,h)=\left\{t\in\mathbb{R}^p\,/\,\left\Vert t-x\right\Vert<h\right\}$, $\Vert \cdot\Vert$ denoting the Euclidean norm of $\mathbb{R}^p$. In what precedes, $\left(k_n\right)_{n\in\mathbb{N}^\ast}$ is a sequence of integers such that $k_n\rightarrow +\infty$ as $n\rightarrow +\infty$. Note that the main difference between this estimator and the usual Parzen-Rosenblatt kernel density estimator is that the bandwith $R_n(x)$ is   random and depends on the $X_i$'s. This estimator is also used for determining and estimator of the regression function $r$. Indeed, assuming that $f(x)>0$, one can easily  see  that 
\begin{equation}\label{r}
r(x)=\frac{g(x)}{f(x)},
\end{equation}
where
\[
g(x)=\int_\mathbb{R} yf_{(X,Y)}(x,y)\mathrm{d}y,
\] 
the function $f_{(X,Y)}$ being the density of $(X,Y)$. Considering the $k$-NN kernel estimator $\widehat{g}_n$ of $g$ defined as 
\begin{equation}\label{gn}
\widehat{g}_{n}(x)=\frac{1}{n(R_{n}(x))^{p}}\sum\limits_{i=1}^{n}Y_iK\left(\frac{X_{i}-x}{R_{n}(x)}\right),\hspace{0.2cm}\textrm{$x\in  \mathbb{R}^{p}$},
\end{equation}
it is seen that by replacing in \eqref{r} $g$ and $f$ by  $\widehat{g}_{n}$ and $\widehat{f}_{n}$ respectively, we obtain an estimator of $r$ which just is the one introduced in \cite{collomb80}. We will modify this estimator as done in \cite{zhu96}, and repeated in  \cite{nkou19}. Specifically,  considering  a sequence $(b_n)_{n\in\mathbb{N}^\ast}$ of positive real numbers converging  to $0$ as $n\rightarrow +\infty$, we define
\[
\widehat{f}_{b_n}(x)=\max( \widehat{f}_n(x),b_n),
\]
and  consider the estimator $\widehat{r}_n$ of  $r$ given by:  
\[
\widehat{r}_{n}(x)=\frac{\widehat{g}_n(x)}{\widehat{f}_{b_n}(x)}.
\]

\section{Rates of uniform consistency}\label{results}
In this section, we present our assumptions, then we give the main results that establish rates of strong
uniform consistency for the estimators of the density and the regression function.

\subsection{Assumptions}
\begin{Assumption}\label{dens}
The density $f$ of $X$ is bounded and bounded from below: there exist $c_0>0$ such that $\inf_{x\in\mathbb{R}^p}f(x)\geqslant c_0$.
\end{Assumption}

\begin{Assumption}\label{lip}
For a given $r\in\mathbb{N}^\ast$, the density $f$ belongs to the class $\mathscr{C}(c,r)$ of functions $\phi\,:\,\mathbb{R}^p\rightarrow\mathbb{R}$ that are $r$ times differentiable and have  $r$-th derivatives $\frac{\partial^r\varphi}{\partial x_{i_1}\cdots \partial x_{i_r}}$, with $(i_1,\cdots,i_r)\in \{1,\cdots,p\}^r$,  satisfying the following Lipschitz condition:
\begin{equation*}
\left\vert\frac{\partial^r\phi}{\partial x_{i_1}\cdots \partial x_{i_r}}(x)-\frac{\partial^r\phi}{\partial x_{i_1}\cdots \partial x_{i_r}}(y)\right\vert\leqslant c\Vert x-y\Vert,
\end{equation*}
where $\Vert\cdot\Vert$ denotes the Euclidean norm of $\mathbb{R}^p$.
\end{Assumption}
\begin{Assumption}\label{reg}
The functions $g_1$ and $g_2$ defined as  $g_1(x)=\mathbb{E}\left( Y \mathds{1}_{\{ Y \geqslant 0 \}} | X = x \right) f(x)$ and   $g_2(x)=\mathbb{E}\left(- Y \mathds{1}_{\{ Y< 0 \}} | X = x \right) f(x)$ are bounded and belong to the class $\mathscr{C}(c,r)$ previously defined.
\end{Assumption}
\begin{Assumption}\label{kernel}
The kernel $K\,:\,\mathbb{R}^p\rightarrow\mathbb{R}$ satisfies the following properties:
\begin{itemize}
\item[(i)]$K$ is bounded, that is $G=\sup_{x\in\mathbb{R}^p}\left\vert K(x)\right\vert<+\infty$.
\item[(ii)]$K$ is symmetric with respect to 0, that is   $K(x)=K(-x)$, $\forall x\in\mathbb{R}^{p}$.
\item[(iii)]$\int_{\mathbb{R}^p}K\left(x\right)\,dx_1\cdots dx_p=1$.
\item[(iv)]$K$ is of order $r$, that is
\[
\int_{\mathbb{R}^p}x_{i_1}\cdots x_{i_\ell}K\left(x\right)\,dx_1\cdots dx_p=0.
\]
for any $\ell\in\{1,\cdots,r\}$ and  $(i_1,\cdots,i_\ell)\in \{1,\cdots,p\}^\ell$.
\item[(v)]
\[
      \int_{\mathbb{R}^{p}}\left\|x\right\|^{r+1}\vert K(x)\vert\,dx_1\cdots dx_p<+\infty.
      \]
\item[(vi)]
$\forall x\in\mathbb{R}^{p}$, $\forall a\in\left[0,1\right]$, $K(ax)\geqslant K(x)$.
\end{itemize}
\end{Assumption}

\begin{Assumption}\label{kn}
 The number $k_{n}$ of neighbors is a sequence of positive integers such that: $k_{n}\sim\lfloor n^{c_{1}}\rfloor$, where $c_{1}\in\left]\frac{1}{2},1\right[$ and $\lfloor a\rfloor$ denotes the integer part of $a$.
\end{Assumption}

\begin{Assumption}\label{bn}
The sequence $(b_n)_{n\in\mathbb{N}^\ast}$ satisfies  $b_{n}\sim n^{-c_{2}}$ with $c_{2}\in\left]0,\frac{1}{10}\right[$.
\end{Assumption}

\begin{Assumption}\label{mn}
There exists a sequence $M_{n}$ of strictly positive numbers such that $M_{n}\sim\sqrt{\log(n)}$ and   $\max_{1\leqslant i\leqslant n}\left\vert Y_{i}\right\vert\leq M_{n}$.
\end{Assumption}

\bigskip

\noindent Assumption \ref{dens} has been made several times in the nonparametric statistics literature. For example, it was made in \cite{zhu07}. Assumption \ref{lip} and \ref{reg} are classical ones; one can found them in  \cite{zhu96}, \cite{zhu07}, \cite{nkou19} for the univariate  case. Assumption \ref{kernel}-$(iv)$ just is   the translation to the multivariate  case of the fact that the kernel $K$ is of   order $r \in\mathbb{N}^\ast$. It is useful in Taylor's expansion used for deriving the consistency rates. Assumption \ref{kernel}-$(vi)$ was made in several works on $k$-NN kernel estimators (e.g., \cite{collomb80,ahmed23}); it is satisfied, for instance, by the Gaussian kernel.  Assumption \ref{mn} is weaker than boundness assumption; for instance, it has  been considered in \cite{nkou23}.

\subsection{Main results}

Now, we give the main results of the paper, that is rates of strong uniform consistency of the estimators introduced in Section \ref{knnest}. First, for the $k$-NN kernel density estimator, we have:
\bigskip

\begin{thm}\label{dens}
Under Assumptions \ref{dens}, \ref{lip}, \ref{kernel} and \ref{kn}, we have:
\begin{equation*}
\sup_{x \in \mathbb{R}^p}\left|\widehat{f}_{n}(x)-f(x)\right| = O_{a.s.}\left(\left(\frac{k_{n}}{n}\right)^{\frac{r+1}{p}}+\sqrt{\frac{n\log(n)}{k_{n}^{2}}}\right).
\end{equation*}
\end{thm}
\bigskip

\begin{remark}
Rates of strong uniform consistency for this estimator was already obtained in \cite{zhang90}, but it was in the univariate case and on any compact subset of $\mathbb{R}$, what is a more restricted framework than ours. In addition, the strong assumption that the kernel has bounded variation on $\mathbb{R}$ was required.
\end{remark}
\bigskip

\noindent In order to get the  rate  for the $k$-NN kernel estimator of the regression function, we first need to deal with the case of the estimator $\widehat{g}_n$ given in \eqref{gn}. We have:

\bigskip

\begin{thm}\label{reg1}
Under Assumptions \ref{dens}, \ref{reg}, \ref{kernel}, \ref{kn} and \ref{mn}, we have:
\begin{equation*}
\sup_{x \in \mathbb{R}^p}\left|\widehat{g}_{n}(x)-g(x)\right| =O_{a.s.}\left(\left(\frac{k_{n}}{n}\right)^{\frac{r+1}{p}}+\sqrt{\frac{n\log(n)M_{n}^{2}}{k_{n}^{2}}}\right). 
\end{equation*}
\end{thm}

\bigskip

\noindent From this result, we obtain as a consequence the following theorem which gives the rate of strong uniform consistency of the $k$-NN kernel estimator of the regression function:

\bigskip

\begin{thm}\label{reg2}
Under Assumptions \ref{dens} to \ref{mn}, we have:
\begin{equation*}
\sup_{x \in \mathbb{R}^p}\left|\widehat{r}_{n}(x)-r(x)\right| =O_{a.s.}\left(\left(\frac{k_{n}}{n}\right)^{\frac{r+1}{p}}+\sqrt{\frac{n\log(n)M_{n}^{2}}{k_{n}^{2}}}+b_n\right).
\end{equation*}
\end{thm}

\section{Proofs of the theorems}\label{preuves}
In this section, we give the proofs of the main results of the paper. First, a result  useful for proving the main theorems is established in Lemma \ref{tech}. Then, the proofs Theorem \ref{dens}, Theorem \ref{reg1} and  Theorem \ref{reg2} are given.
\subsection{A preliminary result}
Let $h\,:\,\mathbb{R}\rightarrow \mathbb{R}$ a measurable function for which there exists a sequence $(\eta_n)_{n\in \mathbb{N}^\ast}$ such that $\vert y\vert\leqslant M_n\Rightarrow \vert h(y)\vert\leqslant \eta_n $, and $1\leqslant\eta_n\leqslant M_n$ for large $n$ enough.    We put
\[
\widehat{\varphi}_n(x)=\frac{1}{n\mathscr{D}_{1,n}^{p}}\sum\limits_{i=1}^{n}h(Y_{i})K\left(\frac{X_i-x}{\mathscr{D}_{2,n}}\right),
\]
where $\left(\mathscr{D}_{1,n} \right)_{n\in \mathbb{N}^\ast}$ and $\left(\mathscr{D}_{2,n} \right)_{n\in \mathbb{N}^\ast}$ are sequences  satisfying
\begin{equation}\label{propdn}
\mathscr{D}_{1,n}^p\geqslant C_1n^{-1}k_n,\,\,\,\,\,\mathscr{D}_{2,n}^p\leqslant C_2n^{-1}k_n,
\end{equation}
for $n$ large enough and some  $C_1>0$ and $C_2>0$, and  
\begin{equation}\label{propdn2}
\left\vert \frac{\mathscr{D}_{2,n}^p}{\mathscr{D}_{1,n}^p}-1\right\vert\sim n^{-\frac{r+1}{p}}.
\end{equation}
Considering
\[
\varphi(x)=\int_\mathbb{R}h(y)f_{(X,Y)}(x,y)\,dy,
\]
where  $f_{(X,Y)}$ is the density of the pair $(X,Y)$, we have:

\bigskip

\begin{lem}\label{tech}
Under the  conditions $(i)$ to $(v)$ of  Assumption \ref{kernel}, if $\varphi$ is bounded and belongs to the class $\mathscr{C}(c,r)$ defined in Assumption \ref{dens}, we have:
\begin{equation*}
\sup_{x\in \mathbb{R}^p}\left\vert\widehat{\varphi}_{n}(x)-\varphi(x)\right\vert=O_{a.s.}\left(\left(\frac{k_{n}}{n}\right)^{\frac{r+1}{p}}+\sqrt{\frac{n\log(n)\eta_n^2}{k_{n}^{2}}}\right).
\end{equation*}
\end{lem}

\noindent\textbf{Proof.}
It suffices to prove the two following properties:
\begin{equation}\label{alea}
\sup_{x\in \mathbb{R}^p}\left\vert\widehat{\varphi}_{n}(x)-\mathbb{E}\left(\widehat{\varphi}_{n}(x) \right)\right\vert =O_{a.s.}\left(\sqrt{\frac{n\log(n)\eta_n^2}{k_{n}^{2}}}\right).
\end{equation}
and
\begin{equation}\label{biais}
\sup_{x\in \mathbb{R}^p}\left\vert\mathbb{E}\left(\widehat{\varphi}_{n}(x)\right)-\varphi(x)\right\vert=O_{a.s.}\left(\left(\frac{k_{n}}{n}\right)^{\frac{r+1}{p}}\right).
\end{equation}

\noindent Proof of \eqref{alea}: From the class of functions
 $$\mathcal{G}_{n}=\left\{\psi_{x}: (t,y)\in\mathbb{R}^p\times [-M_n,M_n]\mapsto\psi_{x}(t,y)=\frac{h(y)}{n\mathcal{D}_{1,n}^p}K\left(\frac{t-x}{\mathcal{D}_{2,n}}\right), x\in \mathbb{R}^p\right\},$$
we use a similar reasoning than in the proof of Theorem 3.1 of \cite{nkou19} (see p. 1299).  Since for any $\psi_x\in\mathcal{G}_{n}$, and for $n$ large enough,
\begin{align*}
 \left|\psi_{x}(t,y)\right|&\leqslant\frac{G\,\vert h(y)\vert}{n\mathcal{D}_{1,n}^p}\leqslant   \frac{G \eta_n}{C_1k_{n}}
\end{align*}
it follows
\[
\mathbb{E}\left[\left|\psi_{x}(X_{j},Y_j)\right|\right] \leqslant\frac{G \eta_n}{C_1k_{n}}=:U_{n},\,\,\,
\mathbb{E}\left[\psi_{x}^{2}(X_{j},Y_j)\right] \leqslant\frac{G ^2\eta_n^2}{C_1^2k_{n}^2}=:\sigma_{n}^{2}.\notag
\]
We can apply Talagrand's inequality (see \cite{talagrand94} and Proposition 2.2 of \cite{gine01}): there exist  $A>0$, $K_{1}>0$ and $K_{2}>0$,  such that for all $t$ satisfying
\begin{align*}
t&\geqslant K_{1}\left[U_{n}\log\frac{AU_{n}}{\sigma_{n}}+\sqrt{n}\sigma_{n}\sqrt{\log\frac{AU_{n}}{\sigma_{n}}}\right]=K_{1}U_{n}\left[\log(A)+\sqrt{n}\sqrt{\log(A)}\right],
\end{align*}
one has
\begin{eqnarray*}
&&P\left\{\sup\limits_{\psi_{x}\in\mathcal{G}_{n}}\left|\sum\limits_{i=1}^{n}\left\{\psi_{x}(X_{i},Y_i)-\mathbb{E}\left(\psi_{x}(X,Y)\right)\right\}\right|> t\right\}\\
&&\leqslant K_{2}\exp\left\{-\frac{1}{K_{2}}\frac{t}{U_{n}}\log\left(1+\frac{tU_{n}}{K_{2}\left(\sqrt{n}\sigma_{n}+U_{n}\sqrt{\log\frac{AU_{n}}{\sigma_{n}}}\right)^{2}}\right)\right\},
\end{eqnarray*}
that is 
\begin{eqnarray}\label{tal}
& &P\left\{\sup\limits_{x\in \mathbb{R}^p}\left|\widehat{\varphi}_{n}(x)-\mathbb{E}\left(\widehat{\varphi}_{n}(x)\right)\right|> t\right\} \nonumber\\
&\leqslant& K_{2}\exp\left\{-\frac{1}{K_{2}}\frac{t}{U_{n}}\log\left(1+\frac{tU_{n}}{K_{2}\left(\sqrt{n}\sigma_{n}+U_{n}\sqrt{\log\frac{AU_{n}}{\sigma_{n}}}\right)^{2}}\right)\right\}\nonumber\\
&= &K_{2}\exp\left\{-\frac{1}{K_{2}}\frac{C_1tk_{n}}{G\eta_n}\log\left(1+\frac{C_1tk_{n}}{nK_{2}G\eta_n\left(1+\sqrt{\log(A)}\right)^2}\right)\right\}.
\end{eqnarray}
Let us put 
$t_{n}=C_3 n^{1/2}k_{n}^{-1} \log^{1/2}(n)\,\eta_n$ and $L=\frac{K_1G}{C_1}\sqrt{\log(A)}$, where
\begin{equation}\label{c1}
 C_{3}> \frac{K_{2}G\left(1+\sqrt{\log(A)}\right)}{C_1}.
\end{equation}
We have,    for $n$ large enough,
\[
C_{3}\log^{1/2}(n)\geqslant 2L\geqslant L\left(1+1\right)\geqslant L\left(1+\sqrt{\frac{\log(A)}{n}}\right) = L\left(\frac{\sqrt{n}+\sqrt{\log(A)}}{\sqrt{n}}\right),
\]
what implies
\[
C_{3}\frac{n^{1/2}\log^{1/2}(n)\eta_n}{k_{n}}\geqslant\frac{K_{1}Gn^{1/2}\eta_n}{C_1k_{n}}\sqrt{\log(A)}\left[\frac{\sqrt{n}+\sqrt{\log(A)}}{\sqrt{n}}\right],
\]
that is $$t_{n}\geqslant K_{1}U_{n}\left[\log(A)+\sqrt{n}\sqrt{\log(A)}\right]= K_{1}\left[U_{n}\log\frac{AU_{n}}{\sigma_{n}}+\sqrt{n}\sigma_{n}\sqrt{\log\frac{AU_{n}}{\sigma_{n}}}\right].$$
Then, \eqref{tal} can be applied to $t_n$ so as to yield $P\left\{\sup\limits_{x\in D}\left|\widehat{\varphi}_{n}(x)-\mathbb{E}\left(\widehat{\varphi}_{n}(x)\right)\right|> t_n\right\} \leqslant u_n$, where
\[
u_n=K_{2}\exp\left\{-\frac{1}{K_{2}}\frac{C_1t_nk_{n}}{G\eta_n}\log\left(1+\frac{C_1t_nk_{n}}{nK_{2}G\eta_n\left(1+\sqrt{\log(A)}\right)^2}\right)\right\}.
\]
Since   $t_nk_n/(n\eta_n)=C_3n^{-1/2}\log^{1/2}(n)\rightarrow 0$ as $n\rightarrow+\infty$, it follows that $u_n\sim v_n$, where 
\begin{align*}
v_n&=K_{2}\exp\left\{-\frac{C_1t_{n}k_{n}}{K_{2}G\eta_n}\times\frac{C_1t_{n}k_{n}}{nK_{2}G\eta_n\left(1+\sqrt{\log(A)}\right)^{2}}\right\}\\
&= K_{2}\exp\left\{-\left(\frac{C_1t_{n}k_{n}n^{-1/2}}{K_{2}G\eta_n\left(1+\sqrt{\log(A)}\right)}\right)^{2}\right\}\\
&= K_{2}\exp\left\{-\left(\frac{C_{1}C_3}{K_{2}G\left(1+\sqrt{\log(A)}\right)}\right)^{2}\log(n)\right\}\\
&=\frac{K_{2}}{n^\alpha},
\end{align*}
with
 $\alpha=\left(\frac{C_{1}C_3}{K_{2}G\left(1+\sqrt{\log(A)}\right)}\right)^{2}.$
From \eqref{c1} we have $\alpha>1$, thus  $\sum_{n=0}^{+\infty}v_n<+\infty$ and $\sum_{n=0}^{+\infty}u_n<+\infty$. Consequently, 
\[
\sum\limits_{n\geqslant 0}P\left\{\sup\limits_{x\in 
\mathbb{R}^p}\left|\widehat{\varphi}_{n}(x)-\mathbb{E}\left(\widehat{\varphi}_{n}(x)\right)\right|>C_{3}\frac{n^{1/2}\log^{1/2}(n)\,\eta_n}{k_{n}}\right\}<+\infty,
\]
and by Borel-Cantelli lemma we  deduce \eqref{alea}.

\bigskip

\noindent Proof of \eqref{biais}: 
\begin{align*}
\mathbb{E}\left(\widehat{\varphi}_{n}(x)\right)&=\frac{1}{\mathscr{D}_{1,n}^{p}}\mathbb{E}\left(h(Y_1)K\left(\frac{X_1-x}{\mathscr{D}_{2,n}}\right)\right)\notag\\
&=\frac{1}{\mathscr{D}_{1,n}^{p}}\int_{\mathbb{R}^{p+1}} h(y)K\left(\frac{t-x}{\mathscr{D}_{2,n}}\right)f_{(X,Y)}(t,y)\,dt_1\cdots dt_p\,dy\\
&=\frac{1}{\mathscr{D}_{1,n}^{p}}\int_{\mathbb{R}^{p}}K\left(\frac{t-x}{\mathscr{D}_{2,n}}\right)\left(\int_\mathbb{R}h(y)f_{(X,Y)}(t,y)\,dy \right) \,dt_1\cdots dt_p\\
&=\frac{1}{\mathscr{D}_{1,n}^{p}}\int_{\mathbb{R}^{p}}K\left(\frac{t-x}{\mathscr{D}_{2,n}}\right)\varphi(t)\,dt_1\cdots dt_p\\
&=\gamma_n\int_{\mathbb{R}^{p}}K\left(u\right)\varphi(x+\mathscr{D}_{2,n}u)\,du_1\cdots du_p,
\end{align*}
where $\gamma_n=\frac{\mathscr{D}_{2,n}^p}{\mathscr{D}_{1,n}^p}$. From  Taylor's theorem, there exists $\theta\in ]0,1[$ such that 
\begin{align*}
\varphi(x+\mathscr{D}_{2,n}u)&=\varphi(x)+\sum\limits_{k=1}^{r-1}\frac{1}{k!}\sum_{1\leqslant i_1,\cdots,i_k\leqslant p}\frac{\partial^k\varphi}{\partial x_{i_1}\cdots \partial x_{i_k}} (x)\,\mathscr{D}_{2,n}^ku_{i_1}\cdots u_{i_k} \\
&\,\,\,\,\,\,\,\,\,\,\,\,\,\,\,\,\,\,\,\,\,\,+\frac{1}{r!}\sum_{1\leqslant i_1,\cdots,i_r\leqslant p}\frac{\partial^r\varphi}{\partial x_{i_1}\cdots \partial x_{i_r}} (x+\theta\mathscr{D}_{2,n}u)\,\mathscr{D}_{2,n}^ru_{i_1}\cdots u_{i_r}.
\end{align*}
Hence
\begin{align*}
&\mathbb{E}\left(\widehat{\varphi}_{n}(x)\right)\\
&=\gamma_n\varphi(x)\int_{\mathbb{R}^{p}}K\left(u\right)\,du_1\cdots du_p\\
&+\gamma_n\sum\limits_{k=1}^{r-1}\frac{1}{k!}\sum_{1\leqslant i_1,\cdots,i_k\leqslant p}\frac{\partial^k\varphi}{\partial x_{i_1}\cdots \partial x_{i_k}} (x)\,\mathscr{D}_{2,n}^k\int_{\mathbb{R}^{p}}u_{i_1}\cdots u_{i_k} K\left(u\right)\,du_1\cdots du_p\\
&+\gamma_n\frac{1}{r!}\sum_{1\leqslant i_1,\cdots,i_r\leqslant p}\mathscr{D}_{2,n}^r\int_{\mathbb{R}^{p}}\frac{\partial^r\varphi}{\partial x_{i_1}\cdots \partial x_{i_r}} (x+\theta\mathscr{D}_{2,n}u)\,u_{i_1}\cdots u_{i_r}\,K\left(u\right)\,du_1\cdots du_p\\
&=\gamma_n\,\varphi(x)\\
&+\gamma_n\frac{1}{r!}\sum_{1\leqslant i_1,\cdots,i_r\leqslant p}\mathscr{D}_{2,n}^r\int_{\mathbb{R}^{p}}\frac{\partial^r\varphi}{\partial x_{i_1}\cdots \partial x_{i_r}} (x+\theta\mathscr{D}_{2,n}u)\,u_{i_1}\cdots u_{i_r}\,K\left(u\right)\,du_1\cdots du_p\\
&=\gamma_n\,\varphi(x)\\
&+\gamma_n\frac{1}{r!}\sum_{1\leqslant i_1,\cdots,i_r\leqslant p}\mathscr{D}_{2,n}^r\int_{\mathbb{R}^{p}}\frac{\partial^r\varphi}{\partial x_{i_1}\cdots \partial x_{i_r}} (x+\theta\mathscr{D}_{2,n}u)\,u_{i_1}\cdots u_{i_r}\,K\left(u\right)\,du_1\cdots du_p\\
&-\gamma_n\frac{1}{r!}\frac{\partial^r\varphi}{\partial x_{i_1}\cdots \partial x_{i_r}} (x)\sum_{1\leqslant i_1,\cdots,i_r\leqslant p}\mathscr{D}_{2,n}^r\int_{\mathbb{R}^{p}}\,u_{i_1}\cdots u_{i_r}\,K\left(u\right)\,du_1\cdots du_p\\
&=\gamma_n\,\varphi(x)+\frac{\gamma_n\mathscr{D}_{2,n}^r}{r!}\\
&\times\sum_{1\leqslant i_1,\cdots,i_r\leqslant p}\int_{\mathbb{R}^{p}}\left(\frac{\partial^r\varphi}{\partial x_{i_1}\cdots \partial x_{i_r}} (x+\theta\mathscr{D}_{2,n}u)-\frac{\partial^r\varphi}{\partial x_{i_1}\cdots \partial x_{i_r}} (x)\right)\,u_{i_1}\cdots u_{i_r}\,K\left(u\right)\,du_1\cdots du_p.
\end{align*}
Since $\varphi$ belongs to $\mathscr{C}(c,r)$, it follows
\begin{align*}
\left\vert\mathbb{E}\left(\widehat{\varphi}_{n}(x)\right)-\varphi(x)\right\vert
&\leqslant\left\vert \gamma_n-1\right\vert\,\,\left\vert\varphi(x)\right\vert\\
& +c\frac{\gamma_n}{r!}\sum_{1\leqslant i_1,\cdots,i_r\leqslant p}\mathscr{D}_{2,n}^{r+1}\theta\int_{\mathbb{R}^{p}}\Vert u\Vert\,\vert u_{i_1}\vert\cdots \vert u_{i_r}\vert\,\vert K\left(u\right)\vert\,du_1\cdots du_p\\
&\leqslant\left\vert \gamma_n-1\right\vert\,\,\left\Vert\varphi\right\Vert_\infty +c\frac{\gamma_n}{r!}p^r\mathscr{D}_{2,n}^{r+1}\int_{\mathbb{R}^{p}}\Vert u\Vert^{r+1}\, \vert K\left(u\right)\vert\,du_1\cdots du_p.
\end{align*}
Since $\gamma_{n}\rightarrow1$ as $n\rightarrow+\infty$, we have for $n$ large enough $\gamma_n\leqslant 3/2$ and, therefore, $\gamma_n\mathscr{D}_{2,n}^{r+1}\leqslant \frac{3}{2}C_2n^{-\frac{r+1}{p}}k_n^{ \frac{r+1}{p}}  $. Thus,
\[
\left\vert\mathbb{E}\left(\widehat{\varphi}_{n}(x)\right)-\varphi(x)\right\vert
\leqslant\left\vert \gamma_n-1\right\vert\,\,\left\Vert\varphi\right\Vert_\infty +C_4\left(\frac{k_{n}}{n}\right)^{\frac{r+1}{p}},
\] 
for some  $C_4>0$. Since $\left\vert \gamma_n-1\right\vert\sim n^{-\frac{r+1}{p}}$ it follows that, for $n$ large enough, $\left\vert\mathbb{E}\left(\widehat{\varphi}_{n}(x)\right)-\varphi(x)\right\vert
\leqslant C_5\left(\frac{k_{n}}{n}\right)^{\frac{r+1}{p}}$ for some  $C_5>0$, what implies \eqref{biais}.
\subsection{Proof of Theorem \ref{dens}}
Considering a sequence $(\beta_n)_{n\in\mathbb{N}^\ast}$ in $]0,1[$ such that $1-\beta_n\sim n^{-\frac{r+1}{p}}$, we put 
\[
D_{n}^{-}(x)=\left[\frac{k_{n}}{nf(x)}\right]^{1/p}\beta_{n}^{1/2p}, \,\,\,D_{n}^{+}(x)=\left[\frac{k_{n}}{nf(x)}\right]^{1/p}\beta_{n}^{-1/2p},
\] 
Then, for $n$ large enough we have almost surely: $D_{n}^{-}(x)\leqslant R_n(x)\leqslant D_{n}^{+}(x)$ (see, e.g., \cite{collomb80}).   According to Assumption \ref{kernel}-$(v)$ we have
\[
K\left(\frac{X_{i}-x}{R_{n}(x)}\right)=K\left(\frac{D_{n}^{-}(x)}{R_{n}(x)}\frac{X_{i}-x}{D_{n}^{-}(x)}\right)\geqslant K\left(\frac{X_{i}-x}{D_{n}^{-}(x)}\right)
\]
and
\[
K\left(\frac{X_{i}-x}{D_{n}^{+}(x)}\right)=K\left(\frac{R_{n}(x)}{D_{n}^{+}(x)}\frac{X_{i}-x}{R_{n}(x)}\right)\geqslant  K\left(\frac{X_{i}-x}{R_{n}(x)}\right).
\]
Thus
\begin{align}\label{ineg}
K\left(\frac{X_{i}-x}{D_{n}^{-}(x)}\right)&\leqslant K\left(\frac{X_{i}-x}{R_{n}(x)}\right)\leqslant K\left(\frac{X_{i}-x}{D_{n}^{+}(x)}\right)
\end{align}
and, therefore, $\widehat{f}_{1,n}(x)\leqslant\widehat{f}_{n}(x)\leqslant\widehat{f}_{2,n}(x)$, where $$\widehat{f}_{1,n}(x)=\frac{1}{n\left(D_{n}^{+}(x)\right)^{p}}\sum\limits_{i=1}^{n}K\left(\frac{X_{i}-x}{D_{n}^{-}(x)}\right)$$
and
$$\widehat{f}_{2,n}(x)=\frac{1}{n\left(D_{n}^{-}(x)\right)^{p}}\sum\limits_{i=1}^{n}K\left(\frac{X_{i}-x}{D_{n}^{+}(x)}\right).$$
Hence
\[
\sup_{x\in \mathbb{R}^p}\left|\widehat{f}_{n}(x)-f(x)\right|\leqslant\max\left\{\sup_{x\in \mathbb{R}^p}\left|\widehat{f}_{1,n}(x)-f(x)\right|,\sup_{x\in \mathbb{R}^p}\left|\widehat{f}_{2,n}(x)-f(x)\right| \right\},
\]
and it remains to prove that
\begin{equation}\label{f1}
\sup_{x\in \mathbb{R}^p}\left|\widehat{f}_{1,n}(x)-f(x)\right| = O_{a.s.}\left(\left(\frac{k_{n}}{n}\right)^{\frac{r+1}{p}}+\sqrt{\frac{n\log(n)}{k_{n}^{2}}}\right) 
\end{equation}
and
\begin{equation}\label{f2}
\sup_{x\in \mathbb{R}^p}\left|\widehat{f}_{2,n}(x)-f(x)\right| = O_{a.s.}\left(\left(\frac{k_{n}}{n}\right)^{\frac{r+1}{p}}+\sqrt{\frac{n\log(n)}{k_{n}^{2}}}\right). 
\end{equation}
For proving \eqref{f1} we apply Lemma \ref{tech} with $h\equiv 1$, $\eta_n\equiv 1$, $\mathscr{D}_{1,n}=D_n^+(x)$ and $\mathscr{D}_{2,n}=D_n^-(x)$. In this case, the properties \eqref{propdn} and \eqref{propdn2} are satisfied. Indeed,  since $\beta_{n}\rightarrow1$ as $n\rightarrow+\infty$, we have for $n$ large enough $1/2\leqslant\beta_n\leqslant 3/2$ and, therefore,
\begin{equation}\label{appl1}
(D_n^+(x))^p\geqslant\sqrt{\frac{2}{3}}\frac{1}{\Vert f\Vert_\infty}n^{-1}k_n,\,\,\,\,(D_n^-(x))^p\leqslant\sqrt{\frac{3}{2}}\frac{1}{c_0}n^{-1}k_n,
\end{equation}
and also
\begin{equation}\label{appl12}
\left\vert \frac{(D_n^-(x))^p}{(D_n^+(x))^p}-1\right\vert=1-\frac{(D_n^-(x))^p}{(D_n^+(x))^p}=1-\beta_n\sim n^{-\frac{r+1}{p}}.
\end{equation}
On the other hand,  $\widehat{\varphi}_n(x)=\widehat{f}_{1,n}(x)$ and
\[
\varphi(x)=\int_\mathbb{R}f_{(X,Y)}(x,y)\,dy=f(x).
\]
Then, applying Lemma \ref{tech} yields \eqref{f1}. Similarly, applying Lemma \ref{tech} to the case where $h\equiv 1$, $\eta_n\equiv 1$, $\mathscr{D}_{1,n}=D_n^-(x)$ and $\mathscr{D}_{2,n}=D_n^+(x)$ leads to \eqref{f2} since
\begin{equation}\label{appl2}
(D_n^-(x))^p\geqslant\frac{1}{\sqrt{2}\Vert f\Vert_\infty}n^{-1}k_n,\,\,\,\,(D_n^+(x))^p\leqslant\frac{\sqrt{2}}{c_0}n^{-1}k_n,
\end{equation}
and 
\begin{equation}\label{appl22}
\left\vert \frac{(D_n^+(x))^p}{(D_n^-(x))^p}-1\right\vert=\frac{(D_n^+(x))^p}{(D_n^-(x))^p}-1=\frac{1-\beta_n}{\beta_n}\sim n^{-\frac{r+1}{p}}.
\end{equation}
\subsection{Proof of Theorem \ref{reg1}}
Clearly, $\widehat{g}_n(x)=\widehat{g}_{1,n}(x)-\widehat{g}_{2,n}(x)$, 
where
\[
\widehat{g}_{1,n}(x)=\frac{1}{n(R_{n}(x))^{p}}\sum\limits_{i=1}^{n}Y_{i}\mathds{1}_{\left\{ Y_{i}\geqslant 0\right\}}K\left(\frac{X_{i}-x}{R_{n}(x)}\right)
\]
and
\[
\widehat{g}_{2,n}(x)=\frac{1}{n(R_{n}(x))^{p}}\sum\limits_{i=1}^{n}(-Y_{i})\mathds{1}_{\left\{Y_{i}<0\right\}}K\left(\frac{X_{i}-x}{R_{n}(x)}\right).
\]
From \eqref{ineg} we get  
$
\widehat{g}_{1,n}^{-}(x)\leqslant\widehat{g}_{1,n}(x)\leqslant \widehat{g}_{1,n}^{+}(x) $
and
$
\widehat{g}_{2,n}^{-}(x)\leqslant\widehat{g}_{2,n}(x)\leqslant \widehat{g}_{2,n}^{+}(x)$, 
where 
$$\widehat{g}_{1,n}^{-}(x)=\frac{1}{n\left(D_{n}^{+}(x)\right)^{p}}\sum\limits_{i=1}^{n}Y_{i}\mathds{1}_{\left\{Y_{i}\geq 0\right\}}K\left(\frac{X_{i}-x}{D_{n}^{-}(x)}\right),$$ $$\widehat{g}_{1,n}^{+}(x)=\frac{1}{n\left(D_{n}^{-}(x)\right)^{p}}\sum\limits_{i=1}^{n}Y_{i}\mathds{1}_{\left\{Y_{i}\geq 0\right\}}K\left(\frac{X_{i}-x}{D_{n}^{+}(x)}\right),$$
$$\widehat{g}_{2,n}^{-}(x)=\frac{1}{n\left(D_{n}^{-}(x)\right)^{p}}\sum\limits_{i=1}^{n}(-Y_{i})\mathds{1}_{\left\{Y_{i}< 0\right\}}K\left(\frac{X_{i}-x}{D_{n}^{+}(x)}\right)$$and $$\widehat{g}_{2,n}^{+}(x)=\frac{1}{n\left(D_{n}^{+}(x)\right)^{p}}\sum\limits_{i=1}^{n}(-Y_{i})\mathds{1}_{\left\{Y_{i}<0\right\}}K\left(\frac{X_{i}-x}{D_{n}^{-}(x)}\right).$$
Then, since $\widehat{g}_n(x)-g(x)=\left(\widehat{g}_{1,n}(x)-g_1(x)\right)-\left(\widehat{g}_{2,n}(x)-g_2(x) \right)$, it follows
\begin{eqnarray*}
& &\sup_{x\in \mathbb{R}^p}\left\vert \widehat{g}_n(x)-g(x)\right\vert\\
& &\leqslant\max\left\{\sup_{x\in \mathbb{R}^p}\left\vert\widehat{g}_{1,n}^{+}(x)-g_1(x)\right\vert+\sup_{x\in \mathbb{R}^p}\left\vert\widehat{g}_{2,n}^{-}(x)-g_2(x)\right\vert,\right.\\
& &\hspace{2cm}\left.\sup_{x\in \mathbb{R}^p}\left\vert\widehat{g}_{1,n}^{-}(x)-g_1(x)\right\vert+\sup_{x\in \mathbb{R}^p}\left\vert\widehat{g}_{2,n}^{+}(x)-g_2(x)\right\vert\right\},
\end{eqnarray*}
and it suffices to show that
\begin{equation}\label{glplus}
\sup_{x\in \mathbb{R}^p}\left\vert\widehat{g}_{\ell,n}^{+}(x)-g_\ell(x)\right\vert=O_{a.s.}\left(\left(\frac{k_{n}}{n}\right)^{\frac{r+1}{p}}+\sqrt{\frac{n\log(n)M_{n}^{2}}{k_{n}^{2}}}\right)
\end{equation}
and
\begin{equation}\label{glmoins}
\sup_{x\in \mathbb{R}^p}\left\vert\widehat{g}_{\ell,n}^{-}(x)-g_\ell(x)\right\vert=O_{a.s.}\left(\left(\frac{k_{n}}{n}\right)^{\frac{r+1}{p}}+\sqrt{\frac{n\log(n)M_{n}^{2}}{k_{n}^{2}}}\right)\hspace{0.2cm}
\end{equation}
for $\ell\in\{1,2\}$. For proving \eqref{glplus} with $\ell=1$,  we apply Lemma \ref{tech} with $h(y)=y\mathds{1}_{\mathbb{R}_+}(y)$, $\eta_n= M_n$, $\mathscr{D}_{1,n}=D_n^-(x)$ and $\mathscr{D}_{2,n}=D_n^+(x)$. In this case, the properties \eqref{propdn} and \eqref{propdn2} are satisfied in \eqref{appl2} and \eqref{appl22} respectively, and we have $\widehat{\varphi}_n(x)=\widehat{g}_{1,n}^{+}(x)$ and
\[
\varphi(x)=\int_\mathbb{R}y\mathds{1}_{\mathbb{R}_+}(y)f_{(X,Y)}(x,y)\,dy=f(x)\int_\mathbb{R}y\mathds{1}_{\mathbb{R}_+}(y)f_{Y\vert X=x}(y)\,dy=g_1(x).
\]
Similarly, applying Lemma \ref{tech} to the case where $h(y)=-y\mathds{1}_{]-\infty,0[}(y)$, $\eta_n= M_n$, $\mathscr{D}_{1,n}=D_n^+(x)$ and $\mathscr{D}_{2,n}=D_n^-(x)$  leads to \eqref{glplus} with $\ell=2$ since  the properties \eqref{propdn} and \eqref{propdn2} are satisfied in \eqref{appl1} and \eqref{appl12} respectively, and we have $\widehat{\varphi}_n(x)=\widehat{g}_{2,n}^{+}(x)$ and
\begin{eqnarray*}
\varphi(x)&=&-\int_\mathbb{R}y\mathds{1}_{]-\infty,0[}(y)f_{(X,Y)}(x,y)\,dy\\
& =&-f(x)\int_\mathbb{R}y\mathds{1}_{]-\infty,0[}(y)f_{Y\vert X=x}(y)\,dy=g_2(x).
\end{eqnarray*}
Equation \eqref{glmoins} is obtained from a similar reasoning.

\subsection{Proof of Theorem \ref{reg2}}
Clearly,
\begin{align*}
\left|\widehat{r}_{n}(x)-r(x)\right|&=\left|\frac{\widehat{g}_{n}(x)}{\widehat{f}_{b_{n}}(x)}-\frac{g(x)}{f(x)}\right|=\frac{\left|\widehat{g}_{n}(x)f(x)-\widehat{f}_{b_{n}}(x)g(x)\right|}{\widehat{f}_{b_{n}}(x)f(x)}\\
&\leqslant c_{_{0}}^{-1}\frac{\left|\left(\widehat{g}_{n}(x)-g(x)\right)f(x)+g(x)\left(f(x)-\widehat{f}_{b_{n}}(x)\right)\right|}{\widehat{f}_{b_{n}}(x)},\\
&\leqslant c_{_{0}}^{-1}\frac{\left\|f\right\|_{\infty} \left\vert\widehat{g}_{n}(x)-g(x)\right\vert+\left\|g\right\|_{\infty}\left\vert\widehat{f}_{b_n}(x)-f(x)\right\vert}{\widehat{f}_{b_{n}}(x)}.
\end{align*}
Since $\sup_{x\in \mathbb{R}^p}\left\vert\widehat{f}_{n}(x)-f(x)\right\vert\rightarrow 0$, a.s., as $n\rightarrow +\infty$, we have for $n$ large enough, $\left\vert\widehat{f}_{n}(x)-f(x)\right\vert\leqslant c_0/2$ and, therefore,
\[
c_0\leqslant f(x)\leqslant \left\vert \widehat{f}_{n}(x)-f(x) \right\vert+\widehat{f}_{n}(x)\leqslant \frac{c_0}{2}+\widehat{f}_{n}(x).
\]
Hence, $\widehat{f}_{n}(x)\geqslant  c_0/2$ and, since $\widehat{f}_{b_n}(x)\geqslant \widehat{f}_{n}(x)$ it follows that $\widehat{f}_{b_n}(x)\geqslant  c_0/2$. Thus
\begin{align*}
\left|\widehat{r}_{n}(x)-r(x)\right|
&\leqslant 2 c_{_{0}}^{-2} \bigg(\left\|f\right\|_{\infty} \left\vert\widehat{g}_{n}(x)-g(x)\right\vert+\left\|g\right\|_{\infty}\left\vert\widehat{f}_{b_n}(x)-f(x)\right\vert\bigg).
\end{align*}
On the other hand, since $\widehat{f}_{n}(x)\leqslant\widehat{f}_{b_n}(x)\leqslant \widehat{f}_{n}(x)+b_n$, it follows that $\left\vert\widehat{f}_{b_n}(x)- \widehat{f}_{n}(x)\right\vert\leqslant b_n$ and, therefore,
\[
\left\vert\widehat{f}_{b_n}(x)-f(x)\right\vert\leqslant\left\vert\widehat{f}_{b_n}(x)- \widehat{f}_{n}(x)\right\vert+\left\vert\widehat{f}_{n}(x)-f(x)\right\vert\leqslant b_n+\left\vert\widehat{f}_{n}(x)-f(x)\right\vert.
\]
Consequently,
\begin{align*}
\sup_{x\in \mathbb{R}^p}\left\vert\widehat{r}_{n}(x)-r(x)\right\vert&\leqslant2 c_{_{0}}^{-2} \bigg(\left\|f\right\|_{\infty}\sup_{x\in \mathbb{R}^p} \left\vert\widehat{g}_{n}(x)-g(x)\right\vert+\left\|g\right\|_{\infty}b_n\\
&\hspace{2cm}+\left\|g\right\|_{\infty}\sup_{x\in \mathbb{R}^p}\left\vert\widehat{f}_{n}(x)-f(x)\right\vert\bigg),
\end{align*}
and the proof is completed by using Theorem \ref{dens} and Theorem \ref{reg1}.


\begin{thebibliography}{99}
                    \bibitem{ahmed23}
		M.S. Ahmed,  M. N'diaye, M.K. Attouch and S. Dabo-Niang, 2023. $k$‑nearest neighbors prediction and classifcation for spatial
data. J. Spatial Econometrics  4, 12.
		\bibitem{attouch17}
		M.  Attouch,  A. Laksaci  and   R. Rafaa, 2017.  Estimation locale lin\'eaire de la r\'egression non param\'etrique
fonctionnelle par la m\'ethode des $k$ plus proches voisins. C. R. Acad. Sci. Paris, Ser. I 355, 824--829.
		\bibitem{bhattacharya90}
                  P.K. Bhattacharya and Y.P. Mack, 1990. Multivariate data-driven $k$-NN function estimation. J. Multivariate  Anal.  35, 1--11.
		\bibitem{burba08}
  F. Burba,  F. Ferraty and P. Vieu, 2008. Convergence de l’estimateur \`a  noyau des $k$  plus proches voisins
en r\'egression fonctionnelle non-param\'etrique . C. R. Acad. Sci. Paris, Ser. I  346, 339--342.
		\bibitem{collomb80} G. Collomb, 1980. Estimation de la r\'egression par la methode des $k$ points les plus proches avec noyau: Quelques propri\'et\'es de convergence ponctuelle. Lectures Notes in Math.  821, 159--175.
		\bibitem{gine01} E. Gin\'e and  A. Guillou, 2001. On consistency of Kernel density estimators for randomly censored data:
rates holding uniformly over adaptive intervals. Ann. Inst. Henri Poincar\'e 	 37, 503--522.
                       \bibitem{kudraszow13} N.L. Kudraszow  and P. Vieu, 2013. Uniform consistency of $k$NN regressors for
functional variables. Statist. Probab. Lett. 	 83, 1863--1870.
                     \bibitem{lian11} H. Lian, 2011. Convergence of functional $k$-nearest neighbor regression estimate with functional responses. Electron. J. Stat.  5, 31--40.
		\bibitem{lu98} Z. Lu  and  P. Cheng, 1998. Strong consistency of nearest neighbor kernel regression
estimation for stationary dependent samples. Sci. China Ser. A  41, 918--926.
                   \bibitem{mack79} Y.P. Mack and M. Rosenblatt, 1979. Multivariate $k$-nearest neighbor density estimates. J. Multivariate  Anal.  9, 1--15.
		\bibitem{moore77} D.S. Moore and J.W. Yackel, 1977. Consistency properties of nearest neighbor density function estimators. Ann. Statist.  5, 143--154.
\bibitem{nkou23} E.D.D. Nkou, 2023. Strong consistency of kernel method for sliced average variance estimation. Comm. Statist. Theory Methods 	 52, 7586--7600.
\bibitem{nkou19} E.D.D. Nkou and G.M. Nkiet, 2019. Strong consistency of kernel estimator in a semiparametric
regression model. Statistics  53, 1289-1305.
\bibitem{talagrand94} M.  Talagrand, 1994. Sharper bounds for Gaussian and empirical processes. Ann.   Probab.  22, 28--76.
\bibitem{zhang90} D.X. Zhang, 1990. Rates of strong uniform convergence of nearest neighbor density estimates on any compact set. Appl. Math. Mech. 11, 385--393.
\bibitem{zhu96} L.X. Zhu  and  K.T. Fang, 1996. Asymptotics for Kernel estimate of sliced inverse regression. Ann. Statist. 	 24, 1053--1068.
\bibitem{zhu07} L. Zhu  and  L.X. Zhu, 2007. On kernel method for sliced average variance estimation. J. Multivariate  Anal. 98, 970--991.
\end{thebibliography}
\end{document}